\documentclass{article}

\usepackage{delarray,verbatim,enumerate,a4wide}
\usepackage{amsmath,amsthm,amstext,amsbsy,amssymb,amsfonts,amscd}
\usepackage[all]{xy}
\usepackage[french]{babel}
\usepackage[T1]{fontenc}
\usepackage[utf8]{inputenc}

\usepackage{scalerel}

\usepackage{upgreek,bbm}

\usepackage[small]{titlesec}

\usepackage{color}
\definecolor{vert}{rgb}{0.1,0.4,0.2}
\usepackage[colorlinks=true,linkcolor=blue,citecolor=vert]{hyperref}

\usepackage{calligra}
\DeclareFontShape{T1}{calligra}{m}{n}{<->s*[0.95]callig15}{}
\DeclareMathAlphabet{\mathscr}{T1}{calligra}{m}{n}

\newtheorem{Th}{Théorème}[]
\newtheorem{Lem}[Th]{Lemme}
\newtheorem{Prop}[Th]{Proposition}
\newtheorem{Cor}[Th]{Corollaire}
\newtheorem{Conj}[Th]{Conjecture}
\newtheorem{Sco}[Th]{Scolie}
\newtheorem{Def} [Th]{Définition}

\newtheorem{DProp}[Th]{Définition \& Proposition}

\newtheorem*{Th*}{Théorème}
\newtheorem*{Lem*}{Lemme}
\newtheorem*{Cor*}{Corollaire}
\newtheorem*{Def*}{Définition}

\def\Preuve{\noindent {\it Preuve.~}}

\def\Remarque{\smallskip\noindent {\it Remarque.~}}
\def\Remarques{\smallskip\noindent {\it Remarques.~}}

		\def\QQ{\mathbb Q}	
\def\NN{\mathbb N}	\def\ZZ{\mathbb Z}		
\def\F2{\mathbb{F}_2}	\def\Z2{\mathbb{Z}_2}		
\def\Zl{{\mathbb{Z}_\ell}} 	\def\Ql{{\mathbb{Q}_\ell}}		

 		\def\P{\mathcal  P}		\def\U{\mathcal  U}	\def\F{\mathcal  F}	
  		\def\C{\mathcal  C}		\def\R{\mathcal  R}	\def\X{\mathcal  X}	\def\K{\mathcal K}
 	  	\def\Cl{\mathcal  C\!\ell}	
\def\E{\mathcal  E}		\def\T{\mathcal  T}			\def\D{\mathcal  D}

				\def\x{{\mathfrak x}}				
		\def\l{{\mathfrak l}}

\def\rg{\operatorname{rg}}	
			
\def\Gal{\operatorname{Gal}}			
\def\Ker{\operatorname{Ker}}	\def\Coker{\operatorname{Coker}}			

\newcommand\scale[2]{\vstretch{#1}{\hstretch{#1}{#2}}}

\newcommand\si[1]{\scale{.7}{#1}}	
\newcommand\ph{{\phantom{*}}}
	\newcommand\lc{{\scale{.8}{\rm lc}}}	
\newcommand\gen{{\scale{.8}{\,\rm gen}}}	\newcommand\cen{{\scale{.8}{\,\rm cen}}}
		
\newcommand\s{{\scale{.7}{\mathcal S}}}		
\def\%{{\scale{.8}{\infty}}}				

\makeatletter
\newcommand*\wt[2][0.2ex]{%
        \begingroup
        \mathchoice{\wt@helper{#1}{#2}{\displaystyle}{\textfont}}
                   {\wt@helper{#1}{#2}{\textstyle}{\textfont}}
                   {\wt@helper{#1}{#2}{\scriptstyle}{\scriptfont}}
                   {\wt@helper{#1}{#2}{\scriptscriptstyle}{\scriptscriptfont}}%
        \endgroup
        #2%
}
\newcommand*\wt@helper[4]{%
        \def\currentfont{\the#41}%
        \def\currentskewchar{\char\the\skewchar\currentfont}%
        \setbox\tw@\hbox{\currentfont$#2$\currentskewchar}%
        \dimen@ii\wd\tw@
        \setbox\tw@\hbox{\currentfont$#2${}\currentskewchar}%
        \advance\dimen@ii-\wd\tw@
        \rlap{\raisebox{-#1}{$\m@th#3\kern\dimen@ii\widetilde{\phantom{#2}}$}}%
}
\makeatother

\def\wE{\,\wt[0.1ex]{\!\mathcal E}}		\def\wU{\wt[0.2ex]{\mathcal U}}
	\def\wCl{\wt[0.1ex]{\mathcal C\!\ell}}

\begin{document}

\title{\Large\bf Sur la trivialité de certains modules d'Iwasawa}
\author{ Jean-François {\sc Jaulent} }
\date{}

\maketitle
\bigskip

{\footnotesize
\noindent{\bf Résumé.} Nous étudions la trivialité de certains modules d'Iwasawa classiques en liaison avec la notion de $\ell$-rationalité pour les corps de nombres totalement $\ell$-adiques.}\smallskip

{\footnotesize
\noindent{\bf Abstract.} We discuss the triviality of some classical Iwasawa modules in connection with the notion of $\ell$-rationality for totally $\ell$-adic number fields.}

\tableofcontents

\section*{Introduction et définitions de base}
\addcontentsline{toc}{section}{Introduction et définitions de base}

Plusieurs des conjectures classiques sur les corps de nombres reviennent à postuler, sous certaines hypothèses, la finitude (ou la pseudo-nullité) d'un module d'Iwasawa convenable.\smallskip

La {\em conjecture de Greenberg} (cf. \cite{Grb1,Grb2}) affirme ainsi que, si $K$ est un corps de nombres totalement réel qui vérifie la conjecture de Leopoldt pour le nombre premier $\ell$  (autrement dit qui admet une unique $\Zl$-extension $K_\%=\bigcup_{n\in\NN} K_n$), la limite projective  $\,\C_{K_\%}\!=\varprojlim \,\Cl_{K_n}$ des $\ell$-groupes de classes d'idéaux attachés aux divers étages de la tour $K_\%/K$ est un groupe fini. Et, sous sa forme généralisée, elle affirme que, si $K$ est un corps de nombres arbitraire vérifiant la conjecture de Leopoldt pour le premier $\ell$ , i.e. admettant exactement $d=c_{\si{K}}+1$ $\Zl$-extensions linéairement indépendantes (où $c_{\si{K}}$ désigne le nombre de places complexes de $K$), la même limite projective $\C_{\bar K}$ prise dans le compositum $\bar K$ de ces $\Zl$-extensions est un module pseudo-nul sur l'algèbre d'Iwasawa $\Lambda_d=\Zl[[\Gal(\bar K/K]]\simeq\Zl[[T_1,\dots,T_d]]$.\smallskip

Or, plus généralement encore, on a:

\begin{Lem}\label{Z}
Tout corps de nombres $K$ de degré $r_{\si{K}}+2c_{\si{K}}$ possède une pro-$\ell$-extension abélienne {\em canonique} $Z$ de groupe de Galois $\Gal(Z/K)\simeq \ZZ_\ell^d$, avec $d=c_{\si{K}}+1$.
\end{Lem}

\Preuve Sous la conjecture de Leopoldt, $Z$ est simplement le compositum des $\Zl$-extensions. Mais, indépendamment de toute conjecture, on peut caractériser $Z$ comme suit: soient $K_\%$ la $\Zl$-extension cyclotomique de $K$, puis $\Gamma=\gamma^\Zl$ le groupe procyclique $\Gal(K_\%/K)$ et $M_\%$ la pro-$\ell$-extension abélienne $\ell$-ramifiée maximale de $K_\%$. Soit alors $\T_\%$ le sous-module de $\Lambda$-torsion du groupe de Galois $\X_\%=\Gal(M_\%/K_\%)$. On sait par la théorie d'Iwasawa (cf. e.g. \cite{Was}) que le quotient $\X_\%/\T_\%$ s'injecte avec un indice fini dans un $\Lambda$-module libre de dimension $c_{\si{K}}$. Alors $Z$ est le sous-corps de $M_\%$ fixé par la racine dans $\X_\%$ du sous-module $\T_\%\X_\%^{\gamma-1}$. En d'autres termes, $Z$ est le compositum des $\Zl$-extensions de $K$ contenues dans le sous-corps de $M_\%$ fixé par $\T_\%\X_\%^{\gamma-1}$.
\medskip

On peut ainsi étendre inconditionnellement la conjecture de Greenberg en postulant:

\begin{Conj}[Conjecture de Greenberg étendue]
Soient $K$ un corps de nombres arbitraire ayant $r_{\si{K}}$ places réelles et $c_{\si{K}}$ places complexes, $Z$ sa pro-$\ell$-extension abélienne canonique de groupe de Galois $\Gal(Z/K)\simeq\Zl^{c_{\si{K}}+\si{1}}$ et $\,\C_Z=\varprojlim \,\Cl_L$ la limite projective (pour les applications normes) des $\ell$-groupes de classes d'idéaux attachés aux sous-corps $L$ de $Z$ de degré fini sur $K$. Le groupe $\,\C_Z$ est alors pseudo-nul comme module sur l'algèbre d'Iwasawa $\Zl[[\Gal(Z/K]]\simeq\Lambda_{c_{\si{K}}+\si{1}}$.
\end{Conj}

Le but de cette note est d'étudier la trivialité du module  $\,\C_Z$, lorsque le corps $K$ est {\em totalement $\ell$-adique} au sens de \cite{J47}, i.e. lorsque ses complétés aux places au-dessus de $\ell$ sont tous de degré 1:

\begin{Def}
Étant donné un nombre premier $\ell$, un corps de nombres $K$ est dit totalement $\ell$-adique lorsque ses complétés aux places au-dessus de $\ell$ sont tous de degré 1 sur $\Ql$, autrement dit lorsque la place $\ell$ est complètement décomposée dans l'extension $K/\QQ$.
\end{Def}

L'intérêt de cette restriction est le lemme de ramification suivant, plus ou moins bien connu, qui joue un rôle essentiel dans notre étude en excluant toute ramification abélienne sauvage au-dessus de la $\Zl$-extension cyclotomique d'un tel corps:

\begin{Lem}\label{cd}
Soient $\ell$ un nombre premier impair, $K$ un corps de nombres totalement $\ell$-adique, $K_\%$\! sa $\Zl$-extension cyclotomique et $N_\%$ une pro-$\ell$-extension de $K_\%$. Si $N_\%$ est localement abélienne sur $K$ aux places au-dessus de $\ell$ (par exemple si $N_\%$ est une pro-$\ell$-extension abélienne de $K$), elle est modérément ramifiée sur $K_\%$ (i.e. non-ramifiée aux places au-dessus de $\ell$).
\end{Lem}

\Preuve C'est une conséquence immédiate de la théorie $\ell$-adique du corps de classes (cf. e.g. \cite{J31}, \S 2.1): le groupe de Galois de la pro-$\ell$-extension abélienne maximale de $\Ql$ s'identifie à la limite projective $\R_{\Ql}=\varprojlim \Ql^\times/\Ql^{\times\ell^n}=(1+\ell)^\Zl\ell^\Zl$ et le sous-groupe d'inertie au groupe procyclique $(1+\ell)^\Zl$. Les places au-dessus de $\ell$ étant totalement ramifiées dans la tour $K_\%/K$, elles ne peuvent plus se ramifier dans $N_\%/K_\%$ dès lors que $N_\%$ est localement abélienne sur $K$.\medskip

Les deux résultats principaux de cette étude (Th. \ref{TP1}, dans le cas totalement réel; Th. \ref{TP2}, dans le cas CM) mettent en avant la $\ell$-rationalité. Rappelons ce dont il s'agit:

\begin{Def}
Un corps de nombres $K$ possédant exactement $c_{\si{K}}$ places complexes est dit $\ell$-rationel pour un nombre premier arbitraire $\ell$ lorsque le groupe de Galois $\Gal(M/K)$ de sa pro-$\ell$-extension abélienne $\ell$-ramifiée $\%$-décomposée (i.e. non-ramifiée en dehors de $\ell$ et complètement décomposée aux places à l'infini) maximale $M$ est un $\Zl$-module libre de dimension $c_{\si{K}}+1$.
\end{Def}

La notion de $\ell$-rationalité d'un corps de nombres $K$, déjà rencontrée par I.R. Shafarevich dans \cite{Sha}, a été formellement introduite par A. Movaheddi dans \cite{Mo} et étudiée en collaboration avec T. Nguyen Quang Do  dans \cite{MN} parallèlement à la notion voisine de $\ell$-régularité introduite par G. Gras et l'auteur dans \cite{GJ} suite aux travaux de \cite{Gra86}. Les deux notions coïncident lorsque $K$ contient le sous-corps réel du corps cyclotomique $\QQ[\zeta_\ell]$ et donnent lieu au même théorème de montée (cf. \cite{J26} pour une synthèse des deux points de vue ou \cite{Gra05} pour plus de détails).\smallskip

Le théorème de Chebotarev permet ainsi de construire étage par étage des $\ell$-tours infinies de corps $\ell$-rationnels en imposant à chaque étape la primitivité de la ramification modérée.\medskip

Dans le cas CM, intervient également la notion de corps logarithmiquement principal (pour le premier donné $\ell$). Le $\ell$-groupe des classes logarithmiques $\,\wCl_K$ a été introduit dans \cite{J28} et son calcul effectif, aujourd'hui implanté directement dans {\sc pari}, est exposé dans \cite{BJ} et \cite{J42}. Il se présente comme un analogue du $\ell$-groupe des classes au sens ordinaire et sa finitude est équivalente à la conjecture de Gross-Kuz'min (cf. e.g. \cite{Gra05,J28, J31,J55}).

\Remarque Les corps $K$ qui vérifient simultanément la conjecture de Leopoldt et celle de Gross-Kuz'min  satisfont plus généralement la conjecture cyclotomique exposée dans \cite{J63}.

\newpage
\section{Trivialité dans la tour cyclotomique}

Intéressons-nous pour commencer au module d'Iwasawa classique attaché à la $\Zl$-extension cyclotomique d'un corps de nombres.

\begin{Th}\label{TP1}
Soient $\ell$ un nombre premier impair, $K$ un corps de nombres totalement $\ell$-adique, $K_\%=\bigcup_{n\in\NN} K_n$ sa $\Zl$-extension cyclotomique et $\,\C_{K_\%}\!=\varprojlim \,\Cl_{K_n}$ la limite projective (pour la norme) des $\ell$-groupes de classes d'idéaux attachés aux divers étages de la tour $K_\%/K$.\par
On a alors $\,\C_{K_\%}\!=1$ si et seulement si le corps $K$ est $\ell$-rationnel et totalement réel, auquel cas il  vérifie banalement les conjectures de Greenberg, de Leopoldt et de Gross-Kuz'min.
\end{Th}

\Preuve Conformément au Lemme \ref{cd}, le corps $K$ étant pris totalement $\ell$-adique, sa pro-$\ell$-extension abélienne $\ell$-ramifiée maximale $M$ est non-ramifiée sur $K_\%$.

 L'égalité  $\,\C_{K_\%}\!=1$ entraîne donc $M=K_\%$. Ainsi $K_\%$ est alors l'unique $\Zl$-extension de $K$, de sorte que $K$ est totalement réel et vérifie la conjecture de Leopoldt (donc aussi celle de Gross-Kuz'min). Enfin la trivialité de $Gal(M/K_\%)$ traduit précisément la $\ell$-rationalité de $K$.

Inversement, si $K$ est totalement réel et $\ell$-rationnel, on a $M=K_\%$; et, plus généralement, $M_n=K_\%$ si $M_n$ désigne la pro-$\ell$-extension abélienne $\ell$-ramifiée de $K_n$, puisque chacun des $K_n$ est encore totalement réel et $\ell$-rationnel en vertu du théorème de propagation de la $\ell$-rationalité donné dans \cite{Gra05,GJ,J26,Mo,MN}. Or, la $\Zl$-extension $K_\%/K$ étant totalement ramifiée (puisque $K$ est pris totalement $\ell$-adique), pour chaque entier $n\in\NN$ le $\ell$-corps de classes de Hilbert $H_n$ de $K_n$ est linéairement disjoint de $K_\%$ sur $K_n$ et on a donc: $\,\Cl_{K_n}\simeq\Gal(H_n/K_n)\simeq\Gal(K_\%H_n/K_\%)=1$, puisque $K_\%H_n$ est contenu dans $M_n$. D'où l'égalité:  $\,\C_{K_\%}\!=1$.\smallskip

\Remarque Dans \cite{Gra21} G. Gras a montré que la condition $\,\C_{K_\%}\!=1$, qui affirme de façon générale la trivialité des $\ell$-goupes de classes $\,\Cl_{K_n}$ des étages finis $K_n$ de la tour cyclotomique $K_\%/K$ pour tout $n$ assez grand, est vérifiée dès lors qu'elle a lieu pour $n=1$, sous réserve que les places au-dessus de $\ell$ soient totalement ramifiées dans la tour, auquel cas elle a lieu pour tout $n\ge1$. En particulier, la condition $\,\C_{K_\%}\!=1$, qui équivaut alors à l'égalité $\,\Cl_{K_1}\!=1$, se lit de ce fait dans $K_1$.\par

Le Théorème \ref{TP1} ci-dessus montre que, sous la condition plus forte de complète décomposition de $\ell$ dans $K/\QQ$, elle se lit directement dans $K$.\medskip

Le théorème de propagation de la $\ell$-rationalité par $\ell$-extension fournit alors un critère nécessaire et suffisant de propagation de la condition de trivialité $\,\C_{K_\%}\!=1$ par $\ell$-extension $\ell$-décomposée:

\begin{Cor}
Soit $L/K$ une $\ell$-extension de corps totalement $\ell$-adiques, $K_\%=\bigcup K_n$ et $L_\%=\bigcup L_n$ leurs $\Zl$-extensions cyclotomiques respectives. Les assertions suivantes sont alors équivalentes:
\begin{itemize}
\item[(i)] Le groupe de Galois $\,\C_{L_\%}=\varprojlim \,\Cl_{L_n}$ est trivial.

\item[(ii)] Le groupe $\,\C_{K_\%}=\varprojlim \,\Cl_{K_n}$ est trivial et l'extension $L/K$ est primitivement ramifiée.
\end{itemize}
\end{Cor}

\Preuve C'est la transposition directe, via le Théorème \ref{TP1}, du théorème de propagation donné dans \cite{GJ,Mo,J26} et \cite{Gra05}. Rappelons qu'une $\ell$-extension $L/K$ est dite {\em primitivement ramifiée} lorsque les logarithmes de Gras (i.e. les images dans le groupe de Galois $\Gal(Z/K)$ du compositum $Z$ des $\Zl$-extensions de $K$) des places {\em modérément} ramifiées dans $L/K$ peuvent être complétées en une $\Zl$-base de $\Gal(Z/K)$.

\begin{Sco}\label{Sco}
Sous les hypothèses du Théorème \ref{TP1}, soit $\,\C'_{K_\%}\!=\varprojlim \,\Cl'_{K_n}$ la limite projective (pour la norme) des $\ell$-groupes de $\ell$-classes des corps $K_n$ (i.e. des quotients respectifs des $\ell$-groupes $\,\Cl_{K_n}$ par leurs sous-groupes engendrés par les classes des idéaux au-dessus de $\ell$).\par
On a alors $\,\C'_{K_\%}\!=1$ si et seulement si le corps $K$ est logarithmiquement principal: $\,\wCl_K=1$.
\end{Sco}

\Preuve En effet, le $\ell$-groupe des classes logarithmiques $\,\wCl_K$ s'interprète par la théorie $\ell$-adique du corps de classes comme groupe de Galois $\Gal(K^\lc/K_\%)$ attaché à la pro-$\ell$-extension abélienne de $K$ localement cyclotomique (i.e. complètement décomposée sur $K_\%$ en chaque place) maximale $K^\lc$. Or, celui-ci  n'est autre que le quotient des genres $^\Gamma\C'_{K_\%}$ de $\,\C'_{K_\%}$ relativement au groupe procyclique $\Gamma=\Gal(K_\%/K)$ (cf. \cite{J28,J31} ou \cite{Gra05}).

\newpage
\section{Trivialité dans le compositum des $\Zl$-extensions}
Soit maintenant $K$ un corps totalement $\ell$-adique de degré $2c_{\si{K}}$, extension quadratique totalement imaginaire d'un sous-corps $K^+$ totalement réel. Notons $M^+$ la pro-$\ell$-extension abélienne $\ell$-ramifiée maximale de $K^+$ et $Z$ la $\ZZ_\ell^{c_{\si{K}}+\si{1}}$-extension canonique de $K$ (cf. Lemme \ref{Z}).\smallskip

Rappelons que, $\ell$ étant impair, si $\bar\tau$ désigne la conjugaison complexe, tout $\Zl[\langle\bar\tau\rangle]$-module $\X$ est somme directe de ses composantes réelle et imaginaire $\X^\pm=\X^{e_\pm}$, avec $e_\pm=\frac{1}{2}(1\pm\bar\tau)$. Ainsi:

\begin{Lem}
Le $\Zl[\langle\bar\tau\rangle]$-module imaginaire $\D_{K_\%}^{\si{[\ell]}-}$ construit sur les idéaux premiers de $K_\%$ au-dessus de $\ell$ est un $\Ql$-espace de dimension $c_{\si{K}}$. Et le sous-module $\P_{K_\%}^{\si{[\ell]}-}$ construit sur les idéaux principaux est un $\Zl$-module libre de même dimension. Il suit:  $\,\Cl_{K_\%}^{\si{[\ell]}-}\simeq(\Ql/\Zl)^{c_{\si{K}}}$.
\end{Lem}

\Preuve Ce résultat est essentiellement bien connu: D'un côté les idéaux premiers au-dessus de $\ell$ étant totalement ramifiés dans la tour cyclotomique, on a: {$\D_{K_n}^{\si{[\ell]}-}=(\D_{K}^{\si{[\ell]}-})^{\ell^{\si{-n}}}$, pour tout $n\in\NN$; d'où, à la limite: $\D_{K_\%}^{\si{[\ell]}-}\simeq\Ql^{c_{\si{K}}}$. D'un autre côté, les $\ell$-unités imaginaires de $K_\%$ provenant directement de $K$, on a, en revanche $\P_{K_n}^{\si{[\ell]}-}=\P_{K}^{\si{[\ell]}-}\simeq\Zl^{c_{\si{K}}}$. D'où: $\,\Cl_{K_\%}^{\si{[\ell]}-}=\D_{K_\%}^{\si{[\ell]}-}/\P_{K_\%}^{\si{[\ell]}-}\simeq(\Ql/\Zl)^{c_{\si{K}}}$.


\begin{Th}\label{TP2}
Soient $\ell$ un nombre premier impair, $K$ un corps de nombres totalement $\ell$-adique, extension quadratique totalement imaginaire d'un sous-corps $K^+$ totalement réel et $Z$ la $\ZZ_\ell^{c_{\si{K}}+\si{1}}$-extension canonique de $K$.
Si le groupe de Galois $\,\C_Z=\Gal(H_Z/Z)$ de la pro-$\ell$-extension abélienne non-ramifiée maximale $H_Z$ de de $Z$ est trivial, on a les conditions suivantes:
\begin{itemize}
\item[(i)] Le corps réel $K^+$ est $\ell$-rationnel et au plus de degré 3.
\item[(ii)] Le $\ell$-groupe des classes logarithmiques de $K$ est trivial: $\,\wCl_K=1$.
\end{itemize}
\end{Th}

\Preuve Supposons $\,\C_Z=1$, introduisons la pro-$\ell$-extension non-ramifiée maximale $\bar H_Z$ de $Z$ et notons $\mathcal G_Z$ son groupe de Galois. Son abélianisé $\,\C_Z$ étant trivial par hypothèse, il en va de même de $\mathcal G_Z$. Ainsi $\bar H_Z$  coïncide avec $Z$. Or, $Z$ étant non-ramifiée sur la $\Zl$-extension cyclotomique $K_\%$ de $K$ en vertu du Lemme \ref{cd}, par construction $\bar H_Z$  est encore la pro-$\ell$-extension non-ramifiée maximale de $K_\%$. En fin de compte $Z$ est donc la réunion $H_\%=\bigcup_{n\in\NN} H_n$ des $\ell$-corps de classes de Hilbert respectifs des étages finis $K_n$ de la tour $K_\%$. Il suit: $\,\C_{K_\%}=\Gal(H_\%/K_\%)=\Gal(Z/H)$.\smallskip

Le  groupe de Galois $\Gal(Z/K_\%)$ étant imaginaire, prenant les composantes réelles, on conclut: $\,\C_{K^+_\%}=\C_{K_\%}^+=1$; et $K^+$ est $\ell$-rationnel en vertu de Théorème \ref{TP1}. 
Il suit: $\,\wCl^+_K=\wCl_{K^+}=1$.\smallskip

Regardons maintenant les composantes imaginaires. Par surjectivité de la norme $\,\C_{K_\%}\to\,\Cl_{K_n}$, nous avons $\rg_\ell\,\Cl^-_{K_n}\le c_{\si{K}}$ pour tout $n\ge 1$, donc $\rg_\ell\,\Cl^-_{K_n}= c_{\si{K}}$ pour $n\ge 1$, en vertu du Lemme. 
Il en résulte que $\,\Cl^-_{K_\%}\simeq(\Ql/\Zl)^{c_{\si{K}}}$ est engendré par les classes des idéaux au-dessus de $\ell$.

En particulier le $\ell$-groupe des $\ell$-classes $\,{\Cl'}_{\!\!K_\%}^-$ est trivial; et, comme il n'y a pas de capitulation pour les $\ell$-classes dans la tour puisque les $\ell$-unités imaginaires sont contenues dans $K$, ce résultat vaut à tous les étages finis: $\,{\Cl'}_{\!\!K_n}^-=1$. On conclut: $\,{\C'}_{\!\!K_\%}^-=\varprojlim \,{\Cl'}_{\!\!K_n}^-=1$; puis
: $\,\wCl_K^-=1$.
\smallskip


Intéressons-nous enfin aux groupes des nœuds respectifs $\K_n$ des $\ell$-extensions abéliennes $H_n/K_n$ (cf. DProp. \ref{noeuds} infra). D'après ce qui précède, les $\ell$-corps de classes de Hilbert respectifs $H'_n$ des corps $H_n$ sont tous contenus dans $H_\%=Z$. Ils sont donc abéliens sur $K$ et l'on a identiquement $H'_n=H_n$ pour chaque $n\in\NN$. Il suit de là que les $\ell$-corps des classes centrales $H_{\si{H_n/K_n}}^\cen$ comme les $\ell$-corps des genres  $H_{\si{H_n/K_n}}^\gen$ coïncident avec les $H_n$; de sorte que dans la suite exacte ($i$) de la Proposition \ref{noeuds} le terme de droite $\Gal(H_{\si{H_n/K_n}}^\cen/H_{\si{H_n/K_n}}^\gen)$ est trivial; d'où par ($ii$) l'isomorphisme:\smallskip

\centerline{$E_{K_n}/E_{K_n}\cap N_{\si{H_n/K_n}}(H_n^\times)\simeq\K_n\simeq\Cl_{K_n}\!\wedge\Cl_{K_n}$,}\smallskip

\noindent puisque l'on a ici: $\Gal(H_n/K_n)\simeq \Cl_{K_n}$ via le corps de classes. Prenant alors les limites projectives pour les applications normes dans la tour cyclotomique $K_\%/K$, on obtient tout comme dans \cite{Fu}, Lem. 3.9, un morphisme surjectif de $\overset{_\leftarrow}{\E}\!_{K_\%}=\varprojlim E_{K_n}$ sur $\,\C_{K_\%}\!\!\wedge\C_{K_\%}$ qui se factorise modulo $\overset{_\leftarrow}{\E}\!_{K_\%}^{\,\gamma-\si{1}}$, puisque $\,\C_{K_\%}$ (et donc son carré alterné) est invariant par $\Gamma=\gamma^\Zl=\Gal(K_\%/K)$.\par


Or,  $\overset{_\leftarrow}{\E}\!_{K_\%}/\overset{_\leftarrow}{\E}\!_{K_\%}^{\,\gamma-\si{1}}$  est un $\Zl$-module de rang essentiel $c_{\si{K}}$ (cf. e.g. \cite{Fu} \S3); et  $\,\C_{K_\%}\!\wedge\,\C_{K_\%}$ est $\Zl$-libre de dimension $\frac{1}{2}c_{\si{K}}(c_{\si{K}}-1)$, puisque $\,\C_{K_\%}\simeq\Gal(Z/K_\%)$ est $\Zl$-libre de dimension $c_{\si{K}}$. Il suit:\smallskip

\centerline{$c_{\si{K}}\ge \frac{1}{2}c_{\si{K}}(c_{\si{K}}-1)$, i.e. $[K^+\!:\QQ]=c_{\si{K}} \le 3$, comme annoncé.}\smallskip


\newpage
\section{Exemple des corps quadratiques totalement $\ell$-adiques}

Pour illustrer les résultats précédents, regardons plus attentivement le cas non-trivial le plus simple: celui des corps quadratiques totalement $\ell$-adiques.\smallskip

Partons donc d'un corps quadratique $K=\QQ[\sqrt d]$, notons $G=\{1,\tau\}$ le groupe $\Gal(K/\QQ)$, prenons un nombre premier impair $\ell$ complètement décomposé dans $K$; écrivons $(\ell)=\l\l'$ dans $K$ et, plus généralement $(\ell)=\l_n^\ph\l_n'$ à chaque étage fini $K_n$ de la $\Zl$-extension cyclotomique $K_n$; notons $M$ la pro-$\ell$-extension abélienne $\ell$-ramifiée et $Z$ le compositum des $\Zl$-extensions de $K$, puis $\,\T_K=\Gal(M/Z)$ et $\,\C_Z$ la limite projective des $\ell$-groupes $\,\Cl_L$ pour $L/K$ de degré fini dans $Z/K$.

\begin{Prop} Pour $K$ quadratique réel, les trois assertions suivantes sont équivalentes:
\begin{itemize}
\item[(i)] Le pro-$\ell$-groupe $\,\C_Z=\,\C_{K_\%}$ est trivial: $\,\C_Z=1$.
\item[(ii)] Le corps $K$ est $\ell$-rationnel: $\,\T_K=1$.
\item[(iii)] On a simultanément $\,\Cl_K=1$ et $\,\wCl_K=1$.
\end{itemize}
\end{Prop}

\Preuve Le corps quadratique réel $K$ admettant pour unique $\Zl$-extension $Z=K_\%$, le Théorème \ref{TP1} nous assure l'équivalence des deux premières assertions. Il reste simplement à vérifier que $\,\T_K$ est trivial si et seulement si $\,\Cl_K$ et $\,\wCl_K$ le sont. \par

Or, d'un côté le $\ell$-corps de classes de Hilbert $H$ de $K$ est linéairement disjoint de $K_\%$, puisque $K_\%/K$ est ici totalement ramifiée, de sorte qu'on a: $\,\Cl_K\simeq\Gal(H/K) \simeq\Gal(HK_\%/K_\%)$. De façon semblable, le $\ell$-groupe des classes logarithmiques vérifie: $\,\wCl_K\simeq\Gal(K^\lc/K_\%$, où $K^\lc$ désigne la pro-$\ell$-extension abélienne localement cyclotomique maximale de $K$. Ainsi, comme $HK_\%$ et $K^\lc$ sont toutes deux contenues dans $M$, on a l'implication: $\,\T_K=1\;\Rightarrow\,\Cl_K=1$ et $\,\wCl_K=1$.\par

D'un autre côté, supposons $\,\Cl_K=1$ et notons $\eta$ un générateur de l'idéal $\l/\l'$ (dans le tensorisé $\R_K=\Zl\otimes_\ZZ K^\times$). La théorie $\ell$-adique du corps de classes (cf. \cite{J31}) nous donne l'isomorphisme $\Gal(M/K^\lc)\simeq\,\wU_\ell/\s_\ell(\wE_K)$, où $\,\wU_\ell=\wU_\l\,\wU_{\l'}$ est le groupe des unités logarithmiques semi-locales et $s_\ell(\wE_K)$ l'image semi-locale du groupe des unités logarithmiques globales. Or, $\,\wU_K$ est engendré ici par les images locales $\ell_\l$ et $\ell_{\l'}$ de $\ell$ dans $\,\wU_\l$ et $\,\wU_{\l'}$, donc conjointement par $\ell_\l\ell_{\l'}=s_\ell(\ell)$ et par $\ell_\l/\ell_{\l'}=s_\ell(\eta)$, c'est-à-dire par l'image de $\,\wE_K$.
Il suit $M=K^\lc$ donc $\,\T_K=1$ pour $\,\wCl_K=1$.

\Remarque La condition $\,\Cl_K=1$ est vérifiée par presque tous les $\ell$ pour $K$ fixé. Il est conjecturé dans \cite{J55} que c'est également le cas de la condition $\,\wCl_K=1$ pour $K$ quadratique réel. La Proposition est donc cohérente avec les heuristiques de Gras qui suggèrent que $K$ est $\ell$-rationnel pour presque tout $\ell$, de sorte que $\,\C_{K_\%}$ serait ainsi presque toujours trivial ici. Pour $\ell$ fixé, en revanche, il résulte de  \cite{Gra22}, \S6 qu'il existe une infinité de corps quadratiques réels $K$ avec $\,\wCl_K\ne 1$.

\begin{Prop}\label{quadim}
Pour $K$ quadratique imaginaire, le pro-$\ell$-groupe $\,\C_Z$ est trivial si et seulement si $K$ est $\ell$-logarithmiquement principal:

\centerline{$\C_Z=1 \Leftrightarrow \,\wCl_K=1$.}
\end{Prop}

\Preuve Le Théorème \ref{TP2} donne l'implication: $\,\C_Z=1 \Rightarrow \,\wCl_K=1$. Reste à vérifier la réciproque. Supposons donc $\,\wCl_K=1$. Comme $\,\wCl_K$ est le quotient des genres de la limite projective des $\ell$-groupes de $\ell$-classes $\,\C'_{K_\%}=\varprojlim \,\Cl'_{K_n}$, cette hypothèse entraîne $\,\C'_{K_\%}=1$, donc finalement $\,\Cl'_{K_n}=1$ pour tout $n\in\NN$. En d'autres termes les $\ell$-groupes de classes $\,\Cl_{K_n}$ sont engendrés par les classes des premiers au-dessus de $\ell$. Plus précisément, puisque sa composante unité $\,\Cl_{\QQ_n}$, qui correspond à l'idempotent $\frac{1}{2}(1+\tau)$, est triviale, $\,\Cl_{K_n}$ est engendré par la classe de l'idéal $\l^\ph_n/\l'_n$ (ou, si l'on préfère, par l'image de la classe $[\l_n]$ par l'idempotent $\frac{1}{2}(1-\tau)$) et le $\ell$-corps de classes de Hilbert $H_n$ de $K_n$ est ainsi une $\ell$-extension cyclique de groupe $G_n\simeq\,\Cl_{K_n}$. Notons $H'_n$ le $\ell$-corps de classes de Hilbert de $H_n$. La formule des classes ambiges de Chevalley (cf. \cite{Chv}) appliquée à l'extension $H_n/K_n$ s'écrit $|\Cl_{H_n}^{\,G_n}|=|\Cl_{K_n}|/[H_n:K_n]=1$ et donne $\,\Cl_{H_n}=1$, i.e. $H'_n=H^\ph_n$. Ainsi $H_n$ est la pro-$\ell$-extension non-ramifiée maximale de $K_n$ et $H_\%=\bigcup_{n\in\NN}H_n$ est celle de $K_\%$.
Maintenant, par le Lemme \ref{cd}, le compositum $Z$ des $\Zl$-extensions de $K$ est contenu dans $H_\%$. Et , comme on a $\Gal(H_\%/K_\%)\simeq\varprojlim\,\Cl_{K_n}\simeq\Zl\simeq\Gal(Z/K_\%)$, l'identité des rangs donne l'égalité $Z=H_\%$.
Il suit $\,\C_Z=1$, comme attendu.

\newpage
\section*{Appendice: Groupe des nœuds, genres et classes centrales}
\addcontentsline{toc}{section}{Appendice: Groupe des nœuds, quotient des genres et classes centrales}

Pour la commodité du lecteur, nous rassemblons ci-dessous quelques résultats classiques sur les relations entre groupe des nœuds et théorie des genres. Pour plus de détails, cf. e.g. \cite{J18}, III.2.1.\smallskip

Pour chaque corps de nombres $K$, nous notons $J_K$ le groupe des idèles, $U_K$ le sous-groupe des idèles unités et $C_K=J_L/K^\times$ le groupe des classes d'idèles. Le corps de classes de Hilbert $H_K$ de $K$, i.e. son extension abélienne non-ramifiée $\%$-décomposée maximale, est ainsi associé au groupe d'idèles $U_K K^\times$ ou, si l'on préfère, au groupe de classes d'idèles $U_K K^\times/K^\times$.

\begin{DProp}\label{genres}
Soit $L/K$ une extension arbitraire de corps de nombres. Alors:\smallskip
\begin{itemize}
\item[(i)] Le compositum $LH_K$ de $L$ avec le corps de classes de Hilbert de $K$ est l'extension abélienne non-ramifiée de $L$ associée au sous-groupe du groupe d'idèles de $L$ défini par:\smallskip

\qquad$J^*_{L/K}= \{ \x \in J_L\;|\; N_{\si{L/K}}(\x)\in U_KK^\times\}$.\smallskip

\item[(ii)] Le corps des genres $H_{\si{L/K}}^\gen$ est la plus grande extension non-ramifiée de $L$ qui provient d'une extension abélienne de $K$. Le sous-groupe d'idèles qui lui correspond est ainsi:\smallskip

\qquad$J^\gen_{\si{L/K}}= \{ \x \in J_L\;|\; N_{\si{L/K}}(\x)\in N_{\si{L/K}}(U_L)K^\times\}$.\smallskip

\item[(iii)] Le corps des classes centrales $H_{\si{L/K}}^\cen$ est, lui, l'extension abélienne de $L$ fixée par:\smallskip

\qquad$J^\cen_{\si{L/K}}= \{ \x \in J_L\;|\; N_{\si{L/K}}(\x)\in N_{\si{L/K}}(U_L L^\times)\}={}_NJ_LU_LL^\times$.\smallskip

Lorsque $L/K$ est galoisienne, $H_{\si{L/K}}^\cen$ est la plus grande extension abélienne non-ramifiée $M$ de $L$, galoisienne sur $K$ et telle que $\Gal(M/L)$ soit contenu dans le centre de $\Gal(M/K)$.\smallskip

\item[(iv)] Tous  sont contenus dans le corps de classes de Hilbert $H_L$ de $L$ fixé par $U_L L^\times$.
\end{itemize}
\end{DProp}


\begin{DProp}\label{noeuds}
Le groupe des nœuds de $L/K$ est le quotient du groupe des normes locales modulo les normes globales: $\K_{L/K}=(K^\times\cap N_{L/K}(J_L))/N_{L/K}(L^\times)$.\smallskip
\begin{itemize}
\item[(i)] De façon générale, $\K_{L/K}$ est relié au groupe de Galois $\Gal(H_{\si{L/K}}^\cen/H_{\si{L/K}}^\gen)$ par la suite exacte:\smallskip

\centerline{$1 \to E_K\cap N_{\si{L/K}}(J_L)/E_K\cap N_{\si{L/K}}(L^\times) \to (K^\times\cap N_{\si{L/K}}(J_L))/N_{\si{L/K}}(L^\times) \to \Gal(H_{\si{L/K}}^\cen/H_{\si{L/K}}^\gen)\to 1$.}\smallskip

\item[(ii)] Et pour $L/K$ abélienne, il s'identifie au quotient du carré alterné de $G=\Gal(L/K)$ par l'image des carrés alternés des sous-groupes de décomposition $G_v$ des places de $K$:\smallskip

\centerline{$\K_{L/K} \simeq (G \wedge G)/\sum_v \phi_v(G_v \wedge G_v)$.}
\end{itemize}
\end{DProp}

La suite exacte ($i$) résulte directement des isomorphismes canoniques:
\begin{center}
$\begin{aligned}
\Gal(H_{\si{L/K}}^\cen/H_{\si{L/K}}^\gen) &\simeq {}^{-\si{1}}N(N(U_L)K^\times/{}_NJ_LU_LL^\times \simeq (N(J_L)\cap N(U_L)K^\times)/N(U_LL^\times)\\ &\simeq (K^\times\cap N(J_L))/ (K^\times\cap N(U_LL^\times))  \simeq (K^\times\cap N(J_L))/N(L^\times)(E_K \cap N(U_L)).
\end{aligned}
$\end{center}

Par ailleurs, lorsque l'extension $L/K$ est galoisienne, J. Tate a donné dans \cite{Ta} une interprétation homologique du groupe  $\K_{L/K}$: la suite exacte de cohomologie associée à la suite courte qui définit le groupe des classes d'idèles $1 \to L^\times \to J_L \to C_L \to 1$ fait apparaître la séquence

\centerline{$ \dots \to \hat{H}^{-\si{1}}(G,J_L) \overset{g}{\to} \hat{H}^{-\si{1}}(G,C_L) \to \hat{H}^{\si{0}}(G,L^\times) \overset{f}{\to}  \hat{H}^{\si{0}}(G,J_L) \to \cdots$}\smallskip

\noindent avec $\K_{L/K} = \Ker f \simeq \Coker g$. Or ce dernier groupe s'interprète via les isomorphismes du corps de classes $\hat{H}^{-\si{1}}(G,C_L) \simeq \hat{H}^{-\si{3}}(G,\ZZ) \simeq H_2(G,\ZZ)$ et $\hat{H}^{-\si{1}}(G,J_L) \simeq \bigoplus_v \hat{H}^{-\si{3}}(G_v,\ZZ) \simeq \bigoplus_v H_2(G_v,\ZZ)$, où, pour chaque place non complexe $v$ de $K$, on désigne par $G_v$ le sous-groupe de décomposition de l'une des places de $L$ au-dessus de $v$.
Cela étant, comme observé par Razar \cite{Raz}, lorsque $G$ est abélien, le groupe d'homologie $H_2(G,\ZZ)$ s'identifie au carré alterné $G \wedge G$ de $G$ et les groupes locaux $H_2(G_v,\ZZ)$ aux carrés alternés $G_v \wedge G_v$ des sous-groupes $G_v$. D'où l'isomorphisme ($ii$).


\Remarques Dans l'isomorphisme obtenu, seules les places ramifiées dans $L/K$interviennent effectivement au dénominateur, puisque pour celles non-ramifiées le sous-groupe de décomposition $G_v$ est cyclique, donc de carré alterné trivial.\par
Enfin, si $L/K$ est une extension galoisienne de corps de nombres, $K'$ une extension de $K$ et $L'=K'L$ l'extension composée, alors dans la description du corps de classes la norme idélique $N_{K'/K}$ correspond à la restriction pour les groupes de Galois (cf. Tate \cite{Ta} ou encore Ozaki \cite{Oz}).
 

\newpage
\def\refname{\normalsize{\sc  Références}}

\medskip\noindent
{\small
\begin{tabular}{l}
Institut de Mathématiques de Bordeaux \\
Université de {\sc Bordeaux} \& CNRS \\
351 cours de la libération\\
F-33405 {\sc Talence} Cedex\\
courriel : Jean-Francois.Jaulent@math.u-bordeaux.fr\\
{\footnotesize \url{https://www.math.u-bordeaux.fr/~jjaulent/}}
\end{tabular}
}

 \end{document}